\documentclass{amsart}
\usepackage{epsf, amssymb, amsfonts, amsmath}%, times, mathptm}

\newtheorem{thm}{Theorem}
\newtheorem{lem}[thm]{Lemma}
\newtheorem{prop}[thm]{Proposition}
\theoremstyle{definition}
\newtheorem{defn}{Definition}
\theoremstyle{remark}

\newcommand{\Q}{{\bf Q}}

\newcommand{\M}{{\mathcal M}}
\newcommand{\R}{{\mathcal R}}
\newcommand{\D}{{\mathcal D}}
\newcommand{\A}{{\mathcal A}}
\newcommand{\I}{{\mathcal I}}
\newcommand{\ka}{{\boldsymbol k}}

\title{{Dimension filtration on loops}}
\author{J.Mostovoy \and J.M. P\'erez-Izquierdo}
\address{Instituto de Matem\'aticas, Unidad Cuernavaca,
 Universidad Nacional Aut\'onoma de M\'exico,
A.P. 273-3,  C.P. 62251, Cuernavaca, Morelos, MEXICO}
\email{jacob@matcuer.unam.mx}
\address{Departamento Matem\'aticas y Computaci\'on, Universidad de La
Rioja, Logro\~no, 26004, SPAIN} \email{jm.perez@dmc.unirioja.es}


\thanks{The second author acknowledges support from
BFM2001-3239-C03-02 (MCYT) and  ANGI2001/26 (Plan Riojano de I+D)}

\begin{document}

\begin{abstract}
We show that the graded group associated to the dimension filtration on
a loop acquires the structure of a Sabinin algebra after being tensored with a
field of characteristic zero. The key to the proof is the interpretation of the
primitive operations of Umirbaev and Shestakov in terms of the operations on 
a loop that measure the failure of the associator to be a homomorphism.
\end{abstract}

\maketitle

\section{Introduction.} 

Let $G$ be a group and $$G=G_{1}\triangleright G_{2}\triangleright G_{3}
\triangleright\ldots$$ --- its lower central series. Then the graded group 
$\oplus G_{i}/G_{i+1}\otimes\Q$ is a Lie algebra with the Lie bracket induced by the 
commutator on $G$. Its universal enveloping algebra can be identified with the algebra 
associated to the filtration of the group ring $\Q G$ by the powers of its augmentation
ideal \cite{Quillen}.  

In this note we generalize these facts to arbitrary loops. 
It will be convenient to speak of dimension series rather than the lower central series.
For groups, both series give rise to the same Lie algebra. 
However, while there are several inequivalent ways of defining the 
lower central series for loops, the definition of the dimension series extends to loops 
without any change.

 \iffalse
Lie algebras arise from associative multiplications
in several different, though closely related, contexts. First and foremost,
Lie algebras are tangent spaces to Lie groups. In the theory of Hopf algebras,
the Lie algebra structure appears as the natural algebraic structure on the 
set of primitive elements of any bialgebra. In the theory of nilpotent groups,
the graded vector space associated to the lower central series or to the
dimension series, has a Lie bracket coming from the group commutator. 

There have been various generalizations of Lie algebras to the non-associative 
context. For example, the infinitesimal theory of Moufang loops produces 
%{\em 
Malcev algebras; %}
non-associative multiplications on symmetric spaces give 
%{\em 
Lie triple systems etc. However, for a relatively long time there was no
satisfactory description of the algebraic structure on tangent spaces to 
general smooth loops. It has been suggested that Akivis algebras
might play the role of Lie algebras for general non-associative 
multiplications. Indeed, the tangent space at the origin of any local
smooth loop has the structure of an Akivis algebra. However, as it turned out,
not every finite-dimensional Akivis algebra comes from a (local) smooth loop. 
\fi

Now it has become clear that 
appropriate generalization of Lie algebras to the non-associative context are
the Sabinin algebras
(called ``hyperalgebras'' in \cite{SM-Orange} and \cite{ShU}).
Sabinin algebras were initially introduced in \cite{SM} as
tangent algebras to local loops. Later, Shestakov and Umirbaev proved that in
any (not necessarily associative) bialgebra the set of primitive elements has
the structure of a Sabinin algebra \cite{ShU}. In fact, it was shown in
\cite{Perez} that any Sabinin algebra can be described as the set of primitive
elements of some bialgebra.

\iffalse
In the present paper we show that Sabinin algebras have their natural place
in the nilpotency theory of loops. Namely, we prove that the graded group
associated to the dimension filtration on a loop acquires the structure of a
Sabinin algebra after being tensored with a field of characteristic zero.
This allows us characterize the primitive elements in the algebra, associated
to the filtration of the loop ring by the powers of its augmentation ideal.
\fi

Our constructions shed some light on the nature of the primitive operations
introduced by Shestakov and Umirbaev in \cite{ShU}. It turns out that the
primitive operations in the free non-associative algebra are induced by
the associator deviations (in the sense of \cite{Mostovoy}) in the free loop.

\section{Dimension subloops.}

Let $L$ be a loop, $\R$ --- a commutative unital ring and
$\R L$ --- the loop ring of $L$ over $\R$. Denote by $IL$ the
augmentation ideal, that is, the kernel of the map
$\R L\to \R$
that sends $\sum a_i x_i$ with $a_i\in\R$ and $x_i\in L$ to $\sum a_i$.
The ideal $IL$ (or simply $I$) is spanned over $\R$ by elements of the
form $x-1$ with $x\in L$. Let $I^mL$ be $m$th power of $I$, that is, the
submodule of $\R L$ spanned over $\R$ by products of at least $m$ elements of
$J$ with any arrangement of the brackets.

\begin{lem} Let $u\in I^m$ and $a\in L$. Then
$au$,  $u/a$ and $a\backslash u$ all lie in $I^m$.
\end{lem}

\begin{proof}
First, $au=(a-1)u+u\in I^m$. Next, $u/a\cdot((a-1)+1)\in I^m$ implies
$u/a\in I$ as $u/a\cdot(a-1)$ is in $I$. This, in turn, implies
$u/a\in I^2$ etc. The same argument works for $a\backslash u$.
\end{proof}

\begin{defn}
The $n$th dimension subloop of $L$ over $\R$ is the intersection
\[D_n(L,\R)=L\cap (1+I^n).\]
\end{defn}
We shall sometimes write $D_n$ instead of $D_n(L,\R)$.

The set $D_n$ is indeed a subloop of $L$. Let $a=1+u$, $b=1+v$ with
$u,v\in I^n$. Clearly, $ab$ is in $D_n$. To see that $a/b$ belongs to $D_n$
take $x=a/b-1$. Then $1+u=(1+v)(1+x)$ and $x=(u-v)/b\in I^n$. In the same
manner one shows that $D_n$ is closed with respect to the left division.

It is clear that the dimension subloops are fully invariant subloops of
$L$ since the augmentation ideal is mapped into itself by any endomorphism
of $\R L$.

\section{The main theorem.}

In this section we will assume that $\R=\ka$ is a field of characteristic
zero. We will denote by $\D$ (or $\D L$) the graded vector space
\[\oplus D_n(L,\ka)/D_{n+1}(L,\ka)\otimes\ka\]
and by $\I$ (or $\I L$) --- the graded algebra
$\oplus I^nL/I^{n+1}L$.

For any loop $L$ the graded algebra $\I$ has a comultiplication
$\I\to \I\otimes\I$ that is a non-trivial algebra homomorphism.
Indeed, the loop algebra $\ka L$ has a comultiplication $\delta$
that sends $g$ in $L$ to $g\otimes g$. Under $\delta$ the element
$g-1$ is sent to
\[ (g-1)\otimes 1+1\otimes (g-1)+(g-1)\otimes (g-1)\]
and, hence, there is an induced comultiplication on the algebra
$\oplus I^n/I^{n+1}$.

There is an inclusion map of $\D$ into $\I$ given by sending
$x\in D_nL$ to $x-1\in I^nL$. The image of the class of $g-1$ is primitive
for any $g\in L$.

Here is our main result:
\begin{thm}\label{main}
The image of $\D$ in $\I$ coincides with the subspace of primitive elements
in $\I$.
\end{thm}

The set of the primitive elements of any bialgebra has the structure of a
Sabinin algebra \cite{ShU}. In fact, if the bialgebra in question is
primitively generated, it can be identified with
the universal enveloping algebra
of the Sabinin algebra of its primitive elements \cite{Perez}.
For any loop the algebra $\I$ is primitively generated since it is generated
by the classes of $g-1\in I$. Therefore, Theorem~\ref{main} can be re-stated
as follows:
\begin{thm}\label{main2}
The graded vector space $\D$ is a Sabinin algebra whose universal enveloping
algebra is $\I$.
\end{thm}

For groups, Theorem~\ref{main2} was first proved by Quillen \cite{Quillen}.
Quillen's result involves the Lie algebra associated to the lower central
series rather than the dimension series; however, for groups these are
isomorphic.

\section{The dimension subloops of a free loop.}

Let $x_i\leftrightarrow x_i'$ be a bijection between two sets of variables
$V$ and $V'$. Denote by $\R[[V']]$ the $\R$-algebra of formal power series in
$m$ non-associative variables $x'_1,\ldots, x'_m$. The power series that start
with $1$ are readily seen to form a loop $\R_0 [[V']]$ under multiplication.
\begin{defn}
The {\em Magnus expansion} is the homomorphism of the free loop $F(V)$ on the
generators $x_1,\ldots, x_m$ into $\R_0 [[V']]$ that sends
the generator $x_i$ to the power series $1+x'_i$.
\end{defn}
We denote the Magnus expansion of $x\in F$ by $\M(x)$. It follows from the
definition that
\[\M(x_i\backslash 1)=1-x'_i+{x'_i}^2-x'_i{x'_i}^2+x'_i(x'_i{x'_i}^2)-\ldots\]
and
\[\M(1/x_i)=1-x'_i+{x'_i}^2-{x'_i}^2x'_i+({x'_i}^2x'_i)x'_i-\ldots.\]

The elements of $F(V)$ whose Magnus expansion has no terms of non-zero degree
less than $n$ form a normal subloop of $F(V)$. These subloops, in fact, are
precisely the dimension subloops:

\begin{lem}\label{m}
The Magnus expansion of $x\in F(V)$ begins with terms of degree $n$ if and
only if $x\in D_n(F(V),\R)$ and $x\notin D_{n+1}(F(V),\R)$.
\end{lem}

It is clear that the Magnus expansion of any element of $I^kF(V)$ starts with
terms of degree at least $k$. The converse is established with the help
of the Taylor formula for free loops.

Define $u_i\in I$ by $u_i=x_i-1$. We shall call the $u_i$ {\em monomials of
degree 1}. A {\em monomial of degree $k$} is a product (with any
arrangement of the brackets) of $k$ monomials of degree 1.
\begin{lem}\label{t}
Given a positive integer $n$, any $x\in F(V)$ can be uniquely written as
\[ x=1+\sum_{j\leq n,\alpha}a_{j,\alpha}\mu_{j,\alpha}+r,\]
where $\mu_{j,\alpha}$ are monomials of degree $j$, $a_{j,\alpha}$ are
elements of $\R$ and $r\in I^{n+1}F(V).$
\end{lem}
The existence of such an expansion follows from the the Taylor formula
with $n=1$. For all the generators $x_i$ it is
obviously true.
Assume now that the set of such $x\in F(V)$ that $x-1$ cannot be written as
a linear combination of the $u_i$ modulo $I^2$, is non-empty. Let $w$ be
a reduced word of minimal possible length (number of {\em operations}
used to form the word) representing such $x$.
If $w=ab$ with $a,b$ reduced words then
\[w-1=ab-1=(a-1)(1+(b-1))+(b-1)\equiv (a-1)+(b-1)\ {\rm mod}\ I^2\]
and we come to a contradiction since $a$ and $b$ are of smaller length then
$w$. If $w=a/b$ with $a,b$ reduced words then
\[w-1=a/b-1\equiv a(1-(b-1))-1 \equiv (a-1)-(b-1) \ {\rm mod}\ I^2\]
since
\[a/b-a(1-(b-1))=(a(1-b))(1-b) \in I^2\]
and we have a contradiction again. Similarly $w$ cannot be of the form
$b\backslash a$; hence $w-1$ is a linear combination of the $u_i$ modulo
$I^2$ for all $w$.

The coefficients  $a_{i,\alpha}$ are uniquely defined as the Magnus expansion
of a monomial $\mu$ starts with
$\mu'$. This proves Lemma~\ref{t}. Lemma~\ref{m} follows immediately from
Lemma~\ref{t}.

\section{Primitive operations and associator deviations.}
Let $\A$ be a bialgebra over a field of characteristic zero, with
non-trivial comultiplication. Assume that $A$
is generated by a set $S$ of its primitive elements. Then the space of all
primitive elements of $\A$ is the minimal vector subspace of $\A$ that contains
$S$ and is closed with respect to the commutatators and the primitive
operations $p_{r,s}$ of Umirbaev and Shestakov \cite{ShU}.

Let $u=(\ldots(x_1x_2)\ldots)x_r$ and $v=(\ldots(y_1y_2)\ldots)y_s$
where $r,s\geq 1$ and $x_i,y_i$ and $z$ are in $\A$.
Specifying the products $u$ and $v$ together with the numbers $r$
and $s$ is equivalent to giving the sequences $x_i$ and $y_i$.
The operations
\[p_{r,s}(x_1,\ldots,x_r;y_1,\ldots,y_s;z)=p_{r,s}(u;v;z)\]
are defined by the formula
\[(uv)z-u(vz)=\sum u_{(1)}v_{(1)}\cdot p_{\alpha,\beta}(u_{(2)};v_{(2)};z)\]
where the sum is taken over all decompositions of the sequences
$x_1,\ldots,x_r$ and $y_1,\ldots,y_s$ into complementary subsequences
$u_{(1)}=(\ldots(x_{i_1}x_{i_2})\ldots)x_{i_p}$,
$u_{(2)}=(\ldots(x_{i_{p+1}}x_{i_{p+2}})\ldots)x_{i_{r}}$
and
$v_{(1)}=(\ldots(y_{j_1}y_{j_2})\ldots)y_{j_q}$,
$v_{(2)}=(\ldots(y_{j_{q+1}}y_{j_{q+2}})\ldots)y_{j_{s}}$.
For example, if $r=s=1$ the operation $p_{1,1}(x_1;y_1;z)$ is just the
associator
\[(x_1,y_1,z)=(x_1y_1)z-x_1(y_1z).\] Also,
\[p_{2,1}(x_1,x_2;y_1;z)=(x_1x_2,y_1,z)-x_1(x_2,y_1,z)-x_2(x_1,y_1,z),\]
\[p_{1,2}(x_1;y_1,y_2;z)=(x_1,y_1y_2,z)-y_1(x_2,y_2,z)-y_2(x_1,y_1,z).\]

Let $\R (V')$ be the free non-associative algebra on $m$
generators $x'_1,\ldots,$ $x'_m$. There is a map
\[\widetilde{\M}: F(V)\to \R (V')\]
defined by taking the lowest-degree term of the Magnus expansion.
Our method of proving Theorem~\ref{main} consists in finding operations
in the free loop $F(V)$ which correspond to the operations $p_{r,s}$ in the
free algebra under the above map.

Such operations were introduced in \cite{Mostovoy} under the name of
{\em associator deviations}. Associator deviations (or simply {\em deviations}) 
are functions from $L^{n+3}$ to $L$
where $L$ is an arbitrary loop and $n$ is a non-negative integer called
the {\em level} of the deviation. There exists one deviation
of level zero, namely the loop associator:
\[(a,b,c)=(a(bc))\backslash ((ab)c).\]
In general,  there are  $(n+2)!/2$ associator deviations of level $n$.
Given $n>0$ and an ordered set $\alpha_1, \ldots, \alpha_n$ of
not necessarily distinct integers satisfying $1\leq\alpha_k\leq n+2$,
the deviation $(a_1,\ldots,a_{n+3})_{\alpha_1, \ldots, \alpha_n}$
is defined inductively by
\[(a_1,\ldots,a_{n+3})_{\alpha_1, \ldots, \alpha_n}:=
(A(a_{\alpha_n})A(a_{\alpha_{n}+1}))\backslash
A(a_{\alpha_n}a_{\alpha_{n}+1}),\]
where $A(x)$ stands for
$(a_1,\ldots,a_{\alpha_n-1},x,a_{\alpha_n+2},
\ldots, a_{n+3})_{\alpha_1, \ldots, \alpha_{n-1}}$.
In particular, there are three deviations of level one:

\[(a,b,c,d)_1=((a,c,d)(b,c,d))\backslash (ab,c,d),\]
\[(a,b,c,d)_2=((a,b,d)(a,c,d))\backslash (a,bc,d),\]
\[(a,b,c,d)_3=((a,b,c)(a,b,d))\backslash (a,b,cd).\]

Let us write $P_{m,n}(x_1,\dots,x_m,y_1,\dots,y_n,z)$ for the deviation
\[(x_1,\dots,x_m,y_1,\dots,y_n,z)_{\underbrace{1,\dots,1}_{m-1},
\underbrace{m+1,\dots,m+1}_{n-1}}.\]
\begin{prop}\label{grancosa}
With the above notation
\begin{multline*}
\M(P_{m,n}(x_1,\dots,x_m,y_1,\dots,y_n,z))\\
= 1 + p_{m,n}(x'_1,\dots,x'_m;y'_1,\dots,y'_n;z') + O(n+m+2)
\end{multline*}
where $O(n+m+2)$ contains no terms of degree less than $n+m+2$.
\end{prop}

The proof of Proposition~\ref{grancosa} is given in the next section. Here
we show how  Proposition~\ref{grancosa} implies Theorem~\ref{main}.

First, let us establish  Theorem~\ref{main} for finitely generated free loops.
Consider the filtration on $\ka (V')$ by be subspaces of elements
of degree at least $i$. The graded algebra associated to this filtration is
clearly $\ka (V')$ itself. It follows from Lemma~\ref{t} that the
map $\widetilde{\M}$ induces an isomorphism between the algebra $\I F(V)$ and
the algebra $\ka (V')$.

Now, assume that all primitive elements on $\I F(V)$ of degree less than $k$
are contained in $\D$; this is certainly true for $k=2$. Any primitive element
of degree $k$ is a linear combination of terms of the form
$p_{\alpha,\beta}(u_1,\ldots,u_{\alpha};v_{1},\ldots,v_{\beta};w)$ and $[u,v]$
where $u,v,u_i,v_i,w$ are primitive elements of $\I F(V)$ of
degree smaller than $k$. Without loss of generality we can assume that
these elements belong to
\[\oplus D_i/D_{i+1}\subset\D F(V)\subset\I F(V);\]
this implies that there exist $\hat{u},\hat{v},\hat{u}_i,
\hat{v}_i$ and $\hat{w}$ in the loop $F(V)$ such that
$\widetilde{\M}(\hat{u})=u$ and similarly for $v,u_i,v_i$ and $w$.
Now, by Proposition~\ref{grancosa}
\[\widetilde{\M}(P_{\alpha,\beta}(\hat{u}_1,\ldots,\hat{u}_{\alpha};
\hat{v}_{1},\ldots,\hat{v}_{\beta};\hat{w})=
p_{\alpha,\beta}(u_1,\ldots,u_{\alpha};v_{1},\ldots,v_{\beta};w). \]
It is also easy to see that
\[\widetilde{\M}((\hat{v}\hat{u})\backslash
(\hat{u}\hat{v}))= [u,v]\]
and, hence, any primitive element of degree $k$ also belongs to $\D$.

Now, let $L$ be an arbitrary loop. Any primitive element $u$ of $\I L$ can be
obtained from a finite number of elements (say, $m$) of the
form $g_i-1\in IL/I^2L$ with $g_i\in L$, by applying commutators, the
$p_{\alpha,\beta}$'s and taking linear combinations.
Consider the homomorphism of the free loop $F(V)$ on $m$ generators to $L$ that
sends the generators of $F(V)$ to the $g_i$.
It is clear that $u$ is the image of
a primitive element $w\in \I F(V)$ under the induced map $\I F(V)\to\I L$.
However, since $1+w\in F(V)$ we see that  $1+u\in L$ and therefore $u$ is in
the image of $\D L$ inside $\I L$.

\section{Proof of  Proposition~\ref{grancosa}.}

Given a subset $S'\subseteq V'$, we shall say that a monomial in $\ka [[V']]$
is {\it balanced} (with respect to $S'$) if it contains each element in
$S'$ at least once. For $S\subseteq V$, we shall say that an
element $\phi \in F(S)$ is balanced (with respect to $S$) if
any nonzero term in $\M(\phi) - 1$ is balanced with respect to $S'$.

\begin{lem}\label{balance}
Given a balanced $\phi \in F(S)$, $x\in S$ (so $\phi = \phi(x)$) and
$y \in V\setminus S$ then the expression
\[\phi(x,y)=(\phi(x)\phi(y))\backslash\phi(xy)\in F(S\cup \{ y\})\]
is balanced too.
\end{lem}

\begin{proof}
\[\M(\phi(xy))-1=(\M(\phi(x))-1)+(\M(\phi(y))-1)+M'\]
where $M'$ is an (infinite) sum of balanced monomials. Hence,
\[\M(\phi(x,y))-1=(\M(\phi(x))\M(\phi(y)))\backslash
(M'-(\M(\phi(x))-1)(\M(\phi(y))-1)).\]
It follows that all the monomials contained in $\M(\phi(x,y))-1$ with
non-zero coefficients are balanced. Indeed, all of the monomials in the
lowest-degree term of $\M(\phi(x,y))-1$ are balanced since every monomial in
\[M'-(\M(\phi(x))-1)(\M(\phi(y))-1)\] is. It follows that
all the other terms of $\M(\phi(x,y))-1$ are balanced since they are
expressed via the lowest-degree term with the help of a recurrent
relation.

\end{proof}

Given $f \in \ka[[V']]$, let $L_{S'}(f)$ be the part of $f$ which contains
all the variables in $S'$ with multiplicity one, but no variables in
$V'\setminus S'$. Similarly, given $\phi \in F(V)$, $L_S(\phi)$ will stand
for $L_{S'}(\M(\phi))$.

\begin{lem}\label{proposition1}
Let $S\subseteq \hat{S}\subseteq V$ with $\vert S\vert \geq 2$,
$x \in S$ and $y\in V\setminus \hat{S}$. Given $\phi = \phi(x) \in F(S)$
balanced,  $\phi(x,y)$ as above and $w \in F(\hat{S})$ then
\[ L_{\hat{S}\cup\{ y\}}(\phi(w,y)) = L_{\hat{S}\cup\{ y\}}(\phi(wy)).\]
\end{lem}

\begin{proof}
First we decompose $S = \{ x\} \sqcup S_0$  ($\sqcup$ denotes disjoint
union),
$\M(\phi(x)) =\sum_I A_I(x)$ with $I$ the
multidegree of $A_I(x)$ on $S'_0$ and $\M(\phi(x,y)) = \sum_K A_K(x,y)$.
With this notation we have that
\[ A_M(xy) = \sum_{I+J+K = M} (A_I(x)A_J(y))A_K(x,y).\]
Hence,
\[ A_M(x,y) = A_M(xy) - A_M(x) - A_M(y) -
\sum_{{I+J+K = M},\ {I,J,K \neq M}}
(A_I(x)A_J(y))A_K(x,y)\]
and, since $\phi(x)$ and $\phi(x,y)$ are balanced,
\[ A_{(1,\dots,1)}(x,y) = A_{(1,\dots,1)}(xy) - A_{(1,\dots,1)}(x) -
A_{(1,\dots,1)}(y).\]
Therefore, using that $ L_{\hat{S}'\cup\{ y'\}}(A_{(1,\dots,1)}(w)) = 0$
and $L_{\hat{S}'\cup\{y'\}}(A_{(1,\dots,1)}(y)) = 0$, we get
\begin{multline*}
L_{\hat{S}\cup\{ y\}}(\phi(w,y)) =
L_{\hat{S}'\cup\{ y'\}}(A_{(1,\dots,1)}(w,y))\\
= L_{\hat{S}'\cup\{
y'\}}(A_{(1,\dots,1)}(wy)) = L_{\hat{S}\cup\{ y\}}(\phi(wy))
\end{multline*}
as desired.

\end{proof}

Now we are in the position to prove Proposition~\ref{grancosa}.
It follows from Lemma~\ref{balance} that
$P_{m,n}(x_1,\dots,x_m,y_1,\dots,y_n,z)$ is
balanced with respect to the set
\[S=\{x_1,\dots,x_m,y_1,\dots,y_n,z\}.\]
Therefore it suffices to show that
\[L_S(P_{m,n}(x_1,\dots,x_m,y_1,\dots,y_n,z)) =
p_{m,n}(x'_1,\dots,x'_m,y'_1,\dots,y'_n,z').\]
As a corollary of Lemma~\ref{proposition1} we have
\[ L_S(P_{m,n}(x_1,\dots,x_m,y_1,\dots,y_n,z)) = L_S((x,y,z)).\]
Set $x' =((x'_1x'_2)\cdots )x'_m$ and
$y'=((y'_1y'_2)\cdots)y'_n$ and recall that
$x=((x_1x_2)\cdots )x_m$ and $y = ((y_1y_2)\cdots)y_n$.
Since $(xy)z = x(yz)\cdot (x,y,z)$, we have
\[
L_S((xy)z)) = \sum_{S_1\sqcup S_2 = S} L_{S_1}(x(yz)) L_{S_2}((x,y,z))
\]
with the convention $L_\emptyset(\cdot) = 1$. The right--hand side of
this equality is
\[x'(y'z') + \sum_{S_1 \sqcup S_2 = S\setminus\{z\}} L_{S_1}(xy)
L_{S_2 \sqcup \{z\}}((x,y,z)) \] since $z'$ appears in any
 term of positive degree of $\M((x,y,z))$. Moreover, setting
\[\hat{p}(x'_{i_1},\dots, x'_{i_k},y'_{j_1},\dots,
y'_{j_l},z') = L_{\{x'_{i_1},\dots, x'_{i_k},y'_{j_1},\dots,
y'_{j_l}\}} (\M((x,y,z)))\]
and $\hat{p}(1,\cdot,z') = \hat{p}(\cdot,1,z')  = \hat{p}(1,1,z') = 0$,
we obtain
\begin{multline*}
(x'y')z' - x'(y'z') = \sum_{S_1
\sqcup S_2 =
S\setminus\{z\}} L_{S_1}(xy) L_{S_2 \sqcup \{z\}}((x,y,z)) \\
= \sum (x'_{(1)}y'_{(1)})\hat{p}(x'_{(2)},y'_{(2)},z'),
\end{multline*}
which shows that the operators $\hat{p}$ agree with the primitive
operations of Umirbaev and Shestakov since they satisfy the same recurrence and
initial conditions. Therefore,
\[
L_S((x,y,z)) = \hat{p}(x'_1,\dots,x'_m,y'_1,\dots,y'_n,z') =
p_{m,n}(x'_1,\dots,x'_m,y'_1,\dots,y'_n,z').
\]

\section{Miscellaneous remarks.}

\subsection*{1.}
All the results of this paper can be stated without change for left
loops (binary systems with left division and a two-sided unit) since
the definition of the associator deviations does not involve right division.

\subsection*{2.}
If all elements of a loop $L$ satisfy some identity, then the
 bialgebra $\I$ satisfies a ``linearized'' version of the same identity.
In particular, if $L$ is Moufang loop, $\D L$ is readily seen to be a
Malcev algebra. Similarly, if $L$ is a Bol loop, $\D L$ is Bol algebra.

\subsection*{3.}
Our results are only stated over fields of characteristic zero
since the necessary general results from the theory of non-associative
algebras are only available for such fields. In particular, the notions
corresponding to Lie rings and restricted Lie algebras are not yet
axiomatized. It is clear, however, that the direct sum
$\oplus D_i(L,\R)/D_{i+1}(L,\R)$ always has a rich algebraic structure.
Using the same argument as in the proof of Lemma~\ref{balance} one sees that
\begin{itemize}
\item $[D_p,D_q]\subset D_{p+q}$;
\item $(D_p,D_q,D_r)\subset D_{p+q+r}$;
\item $(D_{p_1},\ldots,D_{p_n})_{\alpha_1,\ldots,\alpha_{n-3}}
\subset D_{p_1+\ldots+p_n}$ for all combinations
${\alpha_1,\ldots,\alpha_{n-3}}$.
\end{itemize}
Therefore, $\oplus D_i(L,\R)/D_{i+1}(L,\R)$ carries multilinear operations
induced by the commutator and the associator deviations,
and all these operations respect the grading.

\subsection*{4.} For groups, the dimension filtration is closely related to
the lower central series. In particular, the Lie algebras over the field of
rational numbers, associated to both filtrations, are isomorphic. This is
no longer true for general loops, at least if one uses Bruck's definition
of lower central series \cite{Bruck}. Instead, the dimension filtration
is related to the commutator-associator filtration \cite{Mostovoy}; this
connection will be discussed elsewhere.


{\small \begin{thebibliography}{99}
%
\bibitem{Bruck} R.\ H.\ Bruck. A survey of binary systems. Springer-Verlag,
Berlin-Goettingen-Heidelberg, 1958.
%
\bibitem{SM-Orange} P.\ Miheev\ and\ L.\ Sabinin,
{\em Quasigroups and differential geometry},
Quasigroups and loops: theory and applications, 357--430, Heldermann,
Berlin, 1990.
%
\bibitem{Mostovoy} J.\ Mostovoy,
{\em On the notion of the lower central series for loops}, preprint,
available at {\tt http://www.matcuer.unam.mx/\~{}jacob}
%
\bibitem{Perez} J.M.\ P\'erez-Izquierdo,
{\em Algebras, hyperalgebras, nonassociative bialgebras and loops},
preprint, available at {\tt http://mathematik.uibk.ac.at/mathematik/jordan/}
%
\bibitem{SM} L.\ Sabinin, P.\ Mikheev,
{\em Infinitesimal theory of local analytic loops. (Russian)}  Dokl.\ Akad.\
Nauk SSSR  {\bf 297} (1987), 801--804;  translation in  Soviet Math.\ Dokl.\
{\bf 36} (1988), 545--548
%
\bibitem{ShU} I.\ Shestakov\ and\ U.\ Umirbaev,
{\em Free Akivis algebras, primitive elements, and hyperalgebras,}
J. Algebra {\bf 250} (2002), no.~2, 533--548.
%

%\bibitem{Gupta} N.\ Gupta, {\it Free group rings},
%Contemp. Math., 66, Amer. Math. Soc., Providence, RI, 1987.
%
%\bibitem{Hartley} B.\ Hartley,
%{\em Topics in the theory of nilpotent groups},
%in: {Group theory: essays for Philip Hall}. Academic Press,
%London, 1984.
%
\bibitem{Quillen} D.\  Quillen,
{\em On the associated graded ring of a
group ring},  J. Algebra {\bf 10} (1968) 411--418.
%
\end{thebibliography}}
\end{document}